\documentclass[12pt]{article}   
\usepackage{latexsym}        
\usepackage{amssymb}         
\usepackage{amscd}           

\usepackage{amssymb,amsthm,amsmath}

\input epsf

\newtheorem{defn0}{Definition}[section]
\newtheorem{prop0}[defn0]{Proposition}
\newtheorem{thm0}[defn0]{Theorem}
\newtheorem{lemma0}[defn0]{Lemma}
\newtheorem{corollary0}[defn0]{Corollary}
\newtheorem{example0}[defn0]{Example}
\newtheorem{remark0}[defn0]{Remark}

\newtheorem{nota0}[defn0]{Notation}

\newcommand{\defref}[1]{Definition~\ref{#1}}
\newcommand{\propref}[1]{Proposition~\ref{#1}}
\newcommand{\thmref}[1]{Theorem~\ref{#1}}
\newcommand{\lemref}[1]{Lemma~\ref{#1}}
\newcommand{\corref}[1]{Corollary~\ref{#1}}
\newcommand{\exref}[1]{Example~\ref{#1}}

\newcommand{\remref}[1]{Remark~\ref{#1}}

\newcommand{\notref}[1]{Notation~\ref{#1}}

\makeindex

\font\Bbb=msbm10
\newcommand{\ZZ}{\mbox{\Bbb Z}}

\newcommand{\PP}{\mbox{\Bbb P}}

\newcommand{\CC}{\mbox{\Bbb C}}

\newcommand{\AAA}{\mbox{\Bbb A}}

\newcommand{\OO}{\mbox{${\cal O}$}}

\newcommand{\WC}{\mbox{${\cal W}$}}

\newcommand{\sss}{\frak s}

\def\rig#1{\smash{ \mathop{\longrightarrow}
    \limits^{#1}}}

\begin{document}

\begin{center}
 {\Large \bf GROUP ACTION ON INSTANTON BUNDLES OVER  $\PP^3$
}
\end{center}

\begin{center}
LAURA COSTA \footnote{Partially supported by DGICYT PB97-0893.}
\hspace{8mm} GIORGIO OTTAVIANI \footnote{Partially supported by
Italian MURST and by Italian GNSAGA-INDAM.} \end{center}

\centerline{March, 2001}

\vspace{3mm}


{\bf Abstract:} Denote by $MI(k)$ the moduli space of $k$-instanton bundles $E$
 of rank $2$ on $\PP^3=\PP(V)$ and by $Z_k(E)$ the scheme of $k$-jumping lines.
  We prove that $[E]\in MI(k)$ is not stable for the action of $SL(V)$ if
  $Z_k(E)\neq\emptyset$. Moreover $\dim Sym(E)\ge 1$ if $length Z_k(E)\ge 2$.
  We prove also that $E$ is special if and only if $Z_k(E)$ is a smooth conic.
 The action of $SL(V)$ on the moduli of special instanton bundles is studied in detail.

\section{Introduction} \indent

\vspace{3mm}

A $k$-instanton bundle on $\PP^3=\PP(V)$ is a stable $2$-bundle $E$
such that  $c_1=0$, $c_2=k$ and $H^1(E(-2))=0$ (see \cite{Har78} and
\cite{BH78}). The instanton bundles which are trivial
 on the fibers of the twistor map $\PP^3\to S^4$ correspond to self dual
  Yang Mills $Sp(1)$-connections on $S^4$, according to the
ADHM-correspondence.
The moduli space of  $k$-instanton bundles on $\PP^3$ will be denoted
by $MI_{\PP^3}(k)$.

Our aim is to study the natural action of $GL(V)$ over  $MI_{\PP^3}(k)$.
 In this sense this paper can be seen as the natural continuation
 of the study begun in \cite{AO99} trying to answer a general problem
  raised by Simpson who asked about the stable points of the $GL(V)$-action
  on the moduli space of bundles over $\PP(V)$. In \cite{CO00} we proved that
$MI_{\PP^3}(k)$ is an affine variety and this fact simplifies a lot the
geometry of the group action. An easy
consequence is that all points of $MI_{\PP^3}(k)$ are semistable according
to Mumford's geometric invariant theory. For $k\ge 3$ the generic
 $[E]\in MI_{\PP^3}(k)$ has trivial stabilizer $Sym(E)$ and hence,
 for $k\ge 3$, among the nonstable points $[E]$ are all the points $[E]$
 such that $Sym(E)$ has positive dimension. In this work, we will
 focus our attention on the connected component $Sym^0(E)$ containing
  the identity.

A classical tool for the study of vector bundles is the scheme of jumping
lines (see \cite{OSS80}). The restriction of a $k$-instanton bundle $E$ on
a line $L$ is  isomorphic to $\OO_L(-p)\oplus \OO_L(p)$
for some $p$ depending on $L$ and satisfying the inequalities
$0\le p\le k$. It is known that $p=0$ is the general value and if $p\ge 1$
then $L$
is called a jumping line of order $p$. A useful criterion that we prove is the
following (see \thmref{zkunst}, \exref{pitag}):

\vspace{3mm}

\noindent
{\bf Theorem A: }{\em  Let $E$ be a $k$-instanton bundle. If $E$ has a
jumping line of maximal order $k$, then $[E]$ is not stable and
the converse is not true.}

\vspace{3mm}

Moreover, we have found an example of a $k$-instanton bundle $E$ with
$Z_k(E)\neq\emptyset$ and $Sym^0(E)=0$ (see \exref{ce}).

\vspace{3mm}

In \cite{ST90} was observed that the action of $SL(2)\simeq SL(U)$
over $V\simeq U\otimes U$ is useful to get the cohomology and
other geometrical properties of special instanton bundles.
Unfortunately  special instanton bundles are not invariant for
such action. Instead, they are invariant for the action of
$SL(2)\simeq SL(U)$ over $V\simeq U\otimes \CC^2$ where $SL(U)$
acts trivially over $\CC^2$. Our starting point was the remark
that the correct action to be considered is $SL(2)\times
SL(2)\simeq SL(U)\times SL(U')\simeq Spin(V)$
acting over $V\simeq U\otimes U'$
for a suitable isomorphism (see \remref{kernel}), hence for a suitable
symmetric form on $V$. In fact the
symmetry group of the kernel bundle considered in \cite{ST90}
is $SL(U)\times SL(U')$ and
for any special instanton bundle $E$ we get that $Sym(E)$ lies
always in such a group (we believe that this fact is true for any
instanton bundle). This fact conduces us to introduce the
definition of fine action (see \defref{def:properly}). In
particular a subgroup $G\subseteq GL(V)$ has a fine action if
there are $2$-dimensional vector spaces  $U$ and $U'$ such that
$G\subseteq GL(U)\times GL(U')$ and $G$ acts according to fixed
isomorphism $V\simeq U\otimes U'$.

\vspace{3mm}

In Section 5 (see \thmref{big:case} and \thmref{twopoints}) we will prove

\vspace{3mm}

\noindent {\bf Theorem B:} {\em Let $E$ be a $k$-instanton bundle.
If length $Z_k(E)\ge 2$ then $Sym^0(E)$ contains $G$ of positive
dimension having a fine action.}

\vspace{3mm}

In addition, by means of \exref{pitag}, we will see that
the converse of Theorem B is not true. On the other hand, we will see

\vspace{3mm}


\noindent
{\bf Theorem C:}
{\em The following are equivalent
\begin{itemize}
\item[(i)] $Z_k(E)$ is a conic.
\item[(ii)] $E$ is a special $k$-instanton bundle.
\item[(iii)] $Sym^0(E)\supseteq SL(U')$  having a fine  action.
\end{itemize} }

\vspace{3mm}

The equivalence between $(i)$ and $(ii)$ is proved in \corref{line:special},
 $(ii)$ implies $(iii)$ is proved in (\ref{ii-iii}) meanwile
 $(iii)$ implies $(ii)$  is proved in \propref{iii-ii}.

After this result, we are led to conjecture that for any instanton bundle $E$,
 $Sym^0(E)$ has a fine action.

 \vspace{3mm}

In Section 4 we  study the group action over the subvariety
  $MI^s_{\PP^3}(k)$
of special $k$-instanton bundles (see \defref{def:special}) and
our main results are  \corref{semistable} and
\thmref{special}. More precisely, it is known that every
special $k$-instanton bundle is defined by $(k+1)$ lines
through the classical Serre correspondence and that there is a
 pencil of lines whose
branch locus can be studied with the classical machinery of binary forms and
which  is
a $GL(V)$-invariant. In \thmref{special} we show that in general there are
 exactly
$\frac{1}{k+1}\binom{2k}{k}$ $GL(V)$-orbits of special
$k$-instanton bundles with the same branch locus. We mainly apply
results of Trautmann (\cite{Tra88}) that we report in a form
suitable for our purposes.
 Although Trautmann did not apply his results to instanton bundles, probably because
the group action approach on moduli has only nowadays a stronger
importance, some of these applications were certainly known to
him.

We want to point out that we have chosen to present the actions in matrix
form, due to the fact that this is suitable in view of
 computer implementing (\cite{Anc96}).
\vspace{3mm}

\section{Generalities about instantons over $\PP^3$}

\vspace{3mm}

The goal of this section is to prove some general results
about $k$-instanton bundles on $\PP^3$ that will be needed later.
 To this end, we start fixing some notation and recalling some known
 facts (See for instance \cite{BT87}).

\begin{nota0}
\rm $\OO(d)=\OO_{\PP^3}(d)$ denotes the invertible sheaf of degree
$d$ on $\PP^3=\PP(V)$. For any coherent sheaf $E$ on $\PP^3$ we
denote $E(d)=E \otimes \OO_{\PP^3}(d)$ and by $h^iE(d)$ the
dimension of $H^i(\PP^3,E(d))=H^iE(d)$.
\end{nota0}

\begin{defn0}
\label{def:instanton} A $k$-instanton bundle on $\PP^3$ is an
algebraic vector bundle $E$ of rank 2 with Chern classes
$c_1(E)=0$ and $c_2(E)=k$, which is stable (i.e. $ h^0E=0$) and satisfies the
vanishing
$h^1E(-2)=0$. We will denote by $MI_{\PP^3}(k)$ the moduli space
of $k$-instanton bundles on $\PP^3$.
\end{defn0}

\vspace{4mm}

Let $E$ be a $k$-instanton bundle.
It is well known (see \cite{BH78}) that there is a vector space
$I$ of dimension $k$ and a symplectic vector space $(W,J)$
of dimension $2k+2$ such that $E$ is the cohomology of a symplectic monad

\refstepcounter{equation}
 \[ (\theequation) \label{definst} \hspace{8mm}
I^*\otimes\OO(-1)\rig{A^t}W\otimes\OO\rig{A}I\otimes\OO(1).\]

That is, $E\simeq KerA/Im A^t$, where the transpose $A^t$ is computed with
respect to the symplectic form $J$. We have isomorphisms
$I\simeq H^1(E(-1))$, $W\simeq H^1(E\otimes\Omega^1)$.

The vector space $W\otimes I\otimes V=
Hom(W\otimes\OO,I\otimes\OO(1))$ contains the subvariety ${\cal Q}$
 given by
morphisms $A$ such that the sequence (\ref{definst}) is a complex,
that is, such that $AJA^t=0$. In that case,  we shall say that $A$ represents
$E$. $GL(I)\times Sp(W)$ acts on ${\cal Q}$
by $(g,s)\cdot A=gAs$. Let ${\cal Q}^0$ be the open subvariety of
${\cal Q}\subseteq W\otimes I\otimes V$ which consists of morphisms
$A\in{\cal Q}$ that are surjective. By \cite{BH78} (see also
\cite{CO00}) all points of ${\cal Q}^0$ are stable for the action of
$GL(I)\times Sp(W)$ and $MI_{\PP^3}(k)$ is isomorphic to the geometric
quotient
${\cal Q}^0/GL(I)\times Sp(W)$ according to Mumford's GIT. In particular any
morphism between
$k$-instanton bundles lifts to a morphism between the corresponding monads.

\begin{defn0}
\label{unstable:plane} Let $E$ be a $k$-instanton bundle on
$\PP^3$. A hyperplane $ H \subseteq  \PP^3$ is an unstable
hyperplane of $E$ if $h^0E_{|H} \neq 0$. The set of unstable
hyperplanes of $E$ has a natural structure of scheme and it will
be denoted by $W(E)\subseteq{\PP^3}^*$.
\end{defn0}

The fact that, for any $k$-instanton bundle $E$ on
 $\PP^3$, $h^0E(1)\leq 2$ (\cite{BT87}) allows to pose the following
 definition

\begin{defn0}
\label{def:special} A $k$-instanton bundle $E$ on $\PP^3$ is
called special if $ h^0E(1)=2$.
\end{defn0}

\vspace{3mm}

We denote by $MI_{\PP^3}^{s}(k)$ the moduli space of special $k$-instanton
bundles. By its definition, it is constructed as a closed subscheme of
 the affine variety $MI_{\PP^3}(k)$ (\cite{CO00}).
In particular,
$MI_{\PP^3}^{s}(k)$ is also affine  and it has dimension $2k+9$ (see e.g.
  \propref{pi1}).
     Notice that any $E \in MI_{\PP^3}(2)$
is special (\cite{BT87}).

\vspace{3mm}

   For any $k$-instanton bundle $E$ we have $dimW(E)
\le 2$. By a theorem of Coanda (\cite{Coa92}) a $k$-instanton
bundle $E$ is special if and only if $dimW(E)=2$ and in this case
$W(E)$ is a smooth quadric.

\vspace{3mm}

Now we are going to give a similar result  for the scheme of $k$-jumping lines.

\begin{defn0}
\label{jumping} Let $E$ be a $k$-instanton bundle on $\PP^3$ and
$1 \leq p \in \ZZ$. A line $l$ on $\PP^3$ is called a $p$-jumping
line of $E$ if $ E_{|l}= \OO_{l}(-p) \oplus  \OO_{l}(p)$ and we
denote by
\[  Z_j(E)= \{ \quad l \in G(\PP^1,\PP^3)\quad |
 \quad E_{|l}= \OO_{l}(-i) \oplus  \OO_{l}(i)
 \quad j \leq i \quad \}  \]
being  $G(\PP^1,\PP^3)$ the Grassmaniann of lines on $\PP^3$.
\end{defn0}

\noindent $Z_j(E)$ has a natural structure of closed subscheme of
$G(\PP^1,\PP^3)$. From the definition it follows that for any
$k$-instanton bundle $E$ on $\PP^3$ there is the following
filtration
\[G(\PP^1,\PP^3)=Z_0(E) \supseteq Z_1(E) \supseteq Z_2(E).... \]
Moreover, by \cite{Ski97}; Proposition 13 and \cite{Rao97},
$Z_{k+1}(E)=\emptyset$. Hence a line in $Z_k(E)$ is called a maximal order
jumping line.

\begin{remark0}
\label{quadric}
\rm $Z_k(E)$ is obtained by cutting the Pl\"ucker quadric
$Q_4\subseteq\PP^5$ with a linear space. In particular
$deg~Z_k(E)\le 2$. Let $(x_0,\ldots, x_3)$ be homogeneous coordinates
 on $\PP(V)$ and let $A=\sum_{i=0}^3A_ix_i\in W\otimes I\otimes V$
be a matrix representing $E$. Then it is well known that
$Z_k(E)$ is obtained by cutting the
Pl\"ucker quadric $p_{01}p_{23}-p_{02}p_{13}+p_{03}p_{12}=0$
with the linear space $\sum_{i<j}p_{ij}A_iJA_j^t=0$.
\end{remark0}

In the following result we relate the maximal order jumping lines
of $E$ with  unstable hyperplanes of $E$.

\begin{lemma0}
\label{plane:line} Let $E$ be a $k$-instanton bundle on $\PP^3$
and $l \in Z_k(E)$. Then any $H \supseteq l$ is an unstable
hyperplane of $E$, i.e., $H \in W(E)$.
\end{lemma0}
\noindent {\bf Proof.} Since $H \supseteq l$ we have the following
exact sequence
\[ 0\rightarrow E(-1)_{|H} \rightarrow E_{|H}\rightarrow
E_{|l}\rightarrow 0.\] Assume that $H$ is not an unstable plane.
Taking cohomology to the above exact sequence we get the injection
$H^0E_{|l} \hookrightarrow H^1E(-1)_{|H}$. Since by Serre's
duality $h^2E(-1)_{|H}=h^0E(-2)_{|H}=0$, this gives us
\[ k+1=h^0E_{|l} \leq h^1E(-1)_{|H}=-\chi(E(-1)_{|H})=k \]
which is
a contradiction. Hence $H$ is an unstable hyperplane of $E$. \qed

\begin{remark0}
\rm In particular, since for any $k$-instanton bundle $E$ on
$\PP^3$ we have $h^0E_{|H}\leq 1$ (\cite{KO99}; Theorem 1.2), we get
immediately  from the above proof that $Z_j(E)= \emptyset$ for any $j \geq
k+1$.
\end{remark0}

\begin{corollary0}
\label{line:special} Let $E$ be a $k$-instanton bundle on $\PP^3$.
We have $dimZ_k(E)\le 1$. $E$ is special if and only if $dimZ_k(E)=1$
and in this case $Z_k(E)$ is a smooth conic in $G(\PP^1,\PP^3)$.
\end{corollary0}

\noindent {\bf Proof.} It follows from \lemref{plane:line} that if
$dimZ_k(E) \geq 1$ then $dim W(E) \geq 2$. Hence, by
\cite{Coa92}, $E$ is a special instanton bundle. \qed

\begin{remark0}
\rm
 Let $E$ be a $k$-instanton bundle on
$\PP^3$. If $E$ has three different $k$-jumping lines then $E$ is
a special $k$-instanton bundle.
\end{remark0}

\vspace{3mm}

We know (\cite{Har78} for $k=2$ and \cite{BT87} in general) that
 if $E\in MI^s_{\PP^3}(k)$ then $Z_k(E)$ is a smooth
conic on
$G(\PP^1,\PP^3)$ and that $W(E)$ is a smooth quadric. In addition,
$Z_k(E)$ corresponds to one of the two rulings of
$W(E)$ and the generic $s\in H^0E(1)$ vanishes on $k$ disjoint lines
 of the other ruling.
Moreover, it is
well-known that each $E \in MI^s_{\PP^3}(k)$ can be identified
with a $g_{k+1}^1$ without base points on $Z_k(E)\simeq\PP^1$ and that
we can choose this line
as the line of all hyperplanes containing the line $\{x_0=x_1=0\}$.
In fact, the dual of a smooth quadric is isomorphic to the quadric itself,
and we can consider $Q_E\subseteq\PP^3$ as the dual of $W(E)$ and then the
tangent hyperplanes to $Q_E$ belong to $W(E)$. The $g_{k+1}^1$
is constructed by cutting with the zero locus of $s\in H^0E(1)$
on a fixed line of the ruling in $W(E)$ corresponding to $Z_k(E)$. We will
 define the  branch locus of this  $g^1_{k+1}$ as the branch
locus of the bundle. It is defined up to $SL(2)$-action.

\vspace{3mm}

\begin{nota0} \rm
 We will denote by $G^k$ the open subset in
$G(\PP^1,\PP^3)$ which consists of
$g_{k+1}^{1}$ without fixed points and we will denote by
 $\WC \subseteq Gr(\PP^2, \PP^5)$ the open subset of planes
  which cut the Pl\"ucker quadric in a smooth conic.
\end{nota0}

The above discussion shows the following

\begin{prop0}
\label{pi1}
(\cite{HN82}, \cite{BT87};Coroll.4.7)
 The morphism
$\pi_1\colon MI_{\PP^3}^{s}(k)\to \WC$
which takes $E$ to the plane spanned by the conic $Z_k(E)$
is a fibration such that
all the fibers are isomorphic to $G^k$.
\end{prop0}

\vspace{3mm}

Let $Y\subseteq\PP(S^2V)\simeq\PP^9$ be the open subvariety of
smooth quadrics in $\PP(V)$. There is another fibration
$\pi_2\colon MI_{\PP^3}^{s}(k)\to Y$ which takes $E$ to the
quadric $W(E)$. All the fibers of $\pi_2$ are isomorphic to the
disjoint union of two copies of $G^k$ and we have a commutative
diagram
\[ \begin{array}{ccc}
 MI_{\PP^3}^{s}(k)&\rig{\pi_2} &Y\\
||&&\downarrow j\\
 MI_{\PP^3}^{s}(k)&\rig{\pi_1} &\WC\\
\end{array}\]
where $j$ is a $2:1$ finite map. In addition, for any $Q\in \WC$
the fiber $j^{-1}(Q)$ is given by the two rulings of
$Q\simeq\PP^1\times\PP^1$.

\vspace{3mm}

In \cite{AO98}; Proposition 1.1, special $k$-instanton
 bundles were characterized in terms of their monads. However that
description, as the one given in \cite{BT87}, is not useful for our
purposes because the monads considered there were not symplectic
and the techniques that we develop in Section 4
 do not apply.
So we have to introduce an alternative description. With this aim we introduce
some notations

\begin{defn0}
A $(k+1)\times (k+1)$-Hankel matrix $H=(\alpha_{ij})=(\alpha_{i+j-2})$ is
one having equal elements along each diagonal line parallel to the
secondary diagonal, i.e
\[    H= \left (  \begin{array}{cccccc}
     \alpha_0 & \alpha_1 & \alpha_2 &\cdots& \alpha_{k-1} & \alpha_{k} \\
     \alpha_1 & \alpha_2 & \cdots & \cdots & \alpha_{k} & \alpha_{k+1} \\
      \alpha_2 & \vdots   &   & & \alpha_{k+1}& \alpha_{k+2} \\
   \vdots &  \vdots   &     &  &  \vdots       & \vdots \\
    \alpha_{k-1}& \alpha_{k}    &   \cdots &  \cdots   &        \vdots &
\vdots \\
    \alpha_{k} & \alpha_{k+1}& \cdots & \cdots &\alpha_{2k-1}& \alpha_{2k} \\
\end{array} \right ). \]
\end{defn0}

Let $X=[x^k, -kx^{k-1}y, \binom{k}{2}x^{k-2}y^2,\ldots , (-1)^k
y^k]$ and $H$ be a $(k+1)\times (k+1)$-Hankel matrix. Following
\cite{Tra88}; 2.4, we denote by $f_H$ the form of degree $2k$
defined by $X \cdot H \cdot X^t$.

\begin{defn0}
 Denote by $\Delta$ the invariant for
$f_H$  defined by the determinant
\[ \left | \begin{array}{cccc}
  \alpha_0 & \alpha_1 & \cdots      & \alpha_{k}  \\
  \alpha_1 & \alpha_2 & \cdots      & \alpha_{k+1}  \\
   \vdots  &  \vdots  &             & \vdots  \\
  \alpha_{k}     & \alpha_{k+1}    & \cdots       & \alpha_{2k} \\
 \end{array} \right | .   \]
\end{defn0}

Let $(x_0,\ldots, x_3)$ be homogeneous coordinates on $\PP(V)$. We set

\[ I_k(x_0, x_1):= \left ( \begin{array}{llll}
 x_0 & x_1 &       &        \\
     & \ddots& \ddots   &       \\
     &     & x_0&x_1 \\
     \end{array} \right ) \mbox{ and } \quad
 \tilde{I_k}(x_0, x_1):= \left ( \begin{array}{llll}
x_1&&&\\
 x_0 & x_1 &       &        \\
     & \ddots& \ddots   &       \\
     &     & x_0&x_1 \\
&&&x_0\\
     \end{array} \right ). \]

The $k\times (k+1)$-matrix $I_k$ and the $(k+2)\times
(k+1)$-matrix $\tilde{I_k}$ are both well known in the theory of
vector bundles. In fact, the matrix $I_k$ defines a surjective
morphism of vector bundles on $\PP^1=\PP(U)$ given by the natural
multiplication $S^kU\otimes\OO\to S^{k-1}U\otimes \OO (1)$, whose
kernel is isomorphic to $\OO (-k)$. More interesting is the fact
that the converse holds. In fact, it is well known (see e.g.
\cite{AO99}) the following

\vspace{3mm}

\begin{prop0}
\label{represented}
Every surjective morphism of vector bundles on $\PP^1$
\[\OO_{\PP^1}^{k+1}\to\OO_{\PP^1}(1)^k\]
is represented, in a suitable system of coordinates $(x_0, x_1)$,
by the matrix $I_k$.
\end{prop0}

\noindent {\bf Proof.}
The proposition is equivalent to the fact that there is a unique
$SL(2)\times SL(k)\times SL(k+1)$-orbit for nondegenerate
matrices in $\PP(\CC^2\otimes\CC^k\otimes\CC^{k+1})$. \qed

\vspace{3mm}

\begin{nota0}
\label{defbala}
Given  $H$ a nondegenerate $(k+1)\times (k+1)$-Hankel matrix, we denote
\refstepcounter{equation}
 \[ (\theequation) \label{def:balanced} \hspace{8mm}
A:=\left[ I_k(x_0, x_1) | I_k(x_2, x_3)\cdot H\right].\]
\end{nota0}

\vspace{3mm}

\begin{lemma0}
\label{balanced} Let $J$ be the standard skew-symmetric
nondegenerate matrix and let $A$ be as in (\ref{def:balanced}). The following
holds

\begin{itemize}
\item[(i)] $AJA^t=0$.

\item[(ii)] Let $E$ be a $k$-instanton bundle on $\PP^3$ defined by $A$.
Then $h^0E(1)=2$ and $E$ is special.
\item[(iii)]Conversely, every special $k$-instanton bundle is the cohomology
bundle of a monad (\ref{def:balanced}) for a suitable
system of coordinates.
\end{itemize}
\end{lemma0}

\noindent {\bf Proof.} $(i)$ is a straightforward verification.
Let $T_k=\left ( \begin{array}{llll}
&&&1\\
 &&.&        \\
&.&& \\
     1&&&        \\
     \end{array} \right )$ and let
\[ \tilde{A}=\left[ \tilde{I_k}(x_0,x_1)\cdot T_{k+1}\cdot H^{-1}|\tilde {I_k}
(x_2, x_3)\cdot T_{k+1}\right].\]
Then $AJ\tilde{A}^t=0$. Since $\tilde{A}$ has $(k+2)$ rows
it follows $h^0E(1)=2$ which proves $(ii)$.  When $H$ moves among the nondegenerate Hankel matrices,
the moduli space of bundles given by (\ref{def:balanced})
is isomorphic to $\PP(S^{2k}U)\setminus\{\Delta=0\}$ (see \corref{hankel:invariant}). In  \thmref{trautmann}
we quote a result of Trautmann which shows that $\PP(S^{2k}U)\setminus
\{\Delta=0\}$
is isomorphic to
$G^k$. Since this is exactly a fiber of the map $\pi_1$
of \propref{pi1}, we obtain $(iii)$. \qed

\vspace{3mm}

\begin{remark0}
\rm
With the notations of  the above proof set
\[   \tilde{ H}= \left (  \begin{array}{ccccc}
     \alpha_0 & \alpha_1  &\cdots& \alpha_{k} & \alpha_{k+1} \\
     \alpha_1 &  & \cdots  & \alpha_{k+1} & \alpha_{k+2} \\
   \vdots &  \vdots   &       &  \vdots       & \vdots \\
       \alpha_{k-1} & \alpha_{k}& \cdots & \cdots & \alpha_{2k} \\
\end{array} \right ). \]
Then,  $\tilde H\cdot T_{k+2}\cdot \tilde A =A$.
\end{remark0}

\vspace{3mm}

We will end with the following result which is  a straightforward
verification and we left it to the reader.

\begin{lemma0}
Let $E$ be a special $k$-instanton bundle defined by the monad
(\ref{def:balanced}).
Then $W(E)$ is given by all the hyperplanes which are tangent to the quadric
$x_0x_3-x_1x_2$ and $Z_k(E)$ is the conic obtained by cutting the
Pl\"ucker quadric with the plane $\{p_{02}=p_{13}=p_{03}+p_{12}=0\}$. \qed
\end{lemma0}

\section{Useful facts about $SL(2)$-actions}

For any $g \in
\sss l(U)$, we denote by $s^{k}g \in \sss l(S^{k}U)$ the image of
$g$ through the Lie algebra representation of $\sss l(U)$ given by
 the $k$-symmetric power. Notice that $\sss l(U)$ acts on  $S^2(S^kU)$ by
\[ \begin{array}{ccc}
    \sss l(U) \times (S^2(S^kU)) & \longrightarrow & (S^2(S^kU)) \\
    (g,A) & \mapsto & (s^kg)^{t} \cdot A + A \cdot s^{k}g.
    \end{array} \]

 \begin{lemma0}
 \label{action:hankel}
 Let $H$ be a $(k+1)\times (k+1)$-Hankel matrix determined by
 $\alpha=(\alpha_0,...,\alpha_{2k})$ and consider $g \in
\sss l(U)$.  Then the following holds
\begin{itemize}
\item[(i)] $H'=(s^kg)^{t}\cdot H + H \cdot s^{k}g$ is a Hankel matrix.
\item[(ii)] If $\alpha'=(\alpha_0',...,\alpha_{2k}')$ are the coefficients
of $H'$ then
$ (s^{2k}g)^{t} \cdot (\alpha)^{t}=(\alpha')^{t}$.
\item[(iii)] The natural action of  $\sss l(U)$  on forms of degree $2k$
takes $f_H$ to $f_{H'}$.
\end{itemize}
 \end{lemma0}

 \noindent {\bf Proof.} It follows from direct computation. \qed

\begin{corollary0}
\label{hankel:invariant} The $(k+1) \times (k+1)$-Hankel matrices
forms an invariant subspace of $(S^2(S^{k}U))$ isomorphic to
$S^{2k}U$.
\end{corollary0}

\noindent {\bf Proof.} It is an immediate consequence
of \lemref{action:hankel}.
\qed

\vspace{3mm}

There is an analog description for the Lie group $SL(U)$.
If $g\in  SL(U)$,
 we denote by $S^kg\in SL(S^kU)$ the image of
$g$ through the group representation of $SL(U)$ given by the
$k$-symmetric power.

As a consequence of the fact that $S^k(\exp{g})=\exp{(s^kg)}$, from
 \lemref{action:hankel} we obtain

\begin{lemma0}
 \label{group:hankel}
 Let $H$ be a $(k+1)\times (k+1)$-Hankel matrix determined by
 $\alpha=(\alpha_0,...,\alpha_{2k})$ and consider $g \in
SL(U)$.  Then the following holds
\begin{itemize}
\item[(i)] $H''=(S^kg)^{t}\cdot H \cdot S^{k}g$ is a Hankel matrix.
\item[(ii)] If $\alpha''=(\alpha_0'',...,\alpha_{2k}'')$ are the coefficients
of $H''$ then $(S^{2k}g)^{t} \cdot (\alpha)^{t}=(\alpha'')^{t}$.
\item[(iii)] The natural action of  $SL(U)$  on forms of degree $2k$
takes $f_H$ to $f_{H''}$. \qed
\end{itemize}
 \end{lemma0}

\vspace{3mm}

The following is a known result whose generalization can be found in
\cite{AO99}.

\begin{prop0}
\label{actschw}
Let $g\in SL(2)$ be  such that
$g\cdot\left(\begin{array}{c}x'_0\\ x'_1\\ \end{array}\right)=
\left(\begin{array}{c}x_0\\ x_1\\ \end{array}\right)$. Then
\[\left( S^{k-1}g\right)^tI_k(x_0,x_1)\left( S^{k}g^{-1}\right)^t=
I_k(x'_0,x'_1).\]
\end{prop0} \qed

\begin{thm0}   (\cite{Tra88}; Prop. 2.2)
\label{trautmann}
The morphism
 $\PP (S^{2k}U)\setminus\{\Delta=0\}\longrightarrow G^k\subseteq
G(\PP^1,\PP(S^{k+1}U))$
defined by $$ \alpha \mapsto L_{\alpha}=\{(f_0,\ldots ,f_{k+1})\in
S^{k+1}U | (f_0,\ldots ,f_{k+1})\cdot \left ( \begin{array}{cccc}
  \alpha_0 & \alpha_1 & \cdots      & \alpha_{k-1}  \\
  \alpha_1 & \alpha_2 & \cdots      & \alpha_{k}  \\
   \vdots  &  \vdots  &             & \vdots  \\
  \alpha_{k+1}     & \alpha_{k+2}    & \cdots       & \alpha_{2k} \\
 \end{array} \right )=0 \}$$
being $\alpha=(\alpha_0,\ldots ,\alpha_{2k})$, is a
$SL(U)$-equivariant isomorphism. \qed
\end{thm0}

\vspace{3mm}

We have a filtration of $SL(U)$-invariant closed subvarieties
$H_{2k}\subseteq\ldots \subseteq H_2\subseteq H_1=\PP(S^{2k}(U))$
where $H_j=\{f\in S^{2k}(U)|f\hbox{ has a root of multiplicity
}\ge j\}$. In particular, $H_{2k}$ is the rational normal curve of
degree $2k$ and $H_2$ is the discriminant hypersurface. The very
basic example of Mumford's GIT (already known to Hilbert) shows
that $H_{k+1}$ consists of the locus where all the
 $SL(U)$-invariant vanish ( i.e. not semistable points) and the
 $SL(U)$-quotient map
is defined on ${\cal P}^0=\PP(S^{2k}(U))\setminus H_{k+1}$.
We note that $H_{k+1}\subseteq\{\Delta=0\}$.
It is well known that the orbit of $f\in {\cal P}^0$ is not closed in
${\cal P}^0$ if
and only if $f\in H_k$ (i.e. not stable points).
In fact, the orbits of points in $H_k$  all contain in their closure
the orbit of $x^ky^k\in H_k$,
which is the only orbit which is two-dimensional (all the other orbits are
three-dimensional). This particular $SL(U)$-orbit corresponds
to a particular $SL(V)$-orbit of special instanton bundles, that
 we will define in a while.

Now consider the analogous filtration $W_{k+1}\subseteq\ldots
\subseteq W_2\subseteq W_1=\PP(S^{k+1}(U))$ where $W_j=\{f\in
S^{k+1}(U)|f\hbox{ has a root of multiplicity }\ge j\}$. In
particular , $W_{k+1}$ is the rational normal curve of degree
$k+1$ and $W_2$ is the discriminant hypersurface of degree $2k$.

\vspace{3mm}

It is well known (and easy to check) that
a pencil $L\in Gr(\PP^1, \PP(S^{k+1}(U))$ has a fixed point
(that is it does not belong to $G^k$) if
and only if $L$ is contained in some osculating space to the
rational normal curve $W_{k+1}$.
Moreover, the branch locus of $L$ can be identified with
$B_L\in\PP(S^{2k}(U))$ being $B_L$  isomorphic to $L\cap W_2$.
It follows the basic fact that
\[ L\cap W_{s+1}\neq\emptyset\Longleftrightarrow B_L\in H_s
\quad\hbox{for }1\le s\le k\]
and we also remark  that
\[ L_{\alpha}\cap W_{k+1}\neq\emptyset\Longleftrightarrow \alpha\in H_k\]
(compare it with  \cite{Tra88};Prop. 2.6).
The map $G^k\to \PP(S^{2k}(U))$ which takes $L$ to $B_L$
extends to a map
\[R\colon Gr(\PP^1, \PP(S^{k+1}(U))\to \PP(S^{2k}(U))\]
which is finite and equivariant
(see \cite{Tra88}\footnote{The only correction
to be done in the diagram at page 41 of \cite{Tra88} is
that the arrow $H$ must be dotted.};
Remark 2.5, and also our proof of \thmref{special}). In particular,
a pencil $L\in G^k$ is not stable if and only if $L\cap W_{k+1}
\neq\emptyset$. Notice once more
that all the orbits of $L$ such that $L\cap W_{k+1}
\neq\emptyset$ contain in the closure
the unique orbit given by the chordal variety to $W_{k+1}$.

\vspace{3mm}

We resume the above discussion in the following

\begin{corollary0}\label{isogk}
All points of $G^{k}$ are
semistable for the action of $SL(U)$ and the only not stable ones
are the pencils with a point of multiplicity $(k+1)$. All the not
stable orbits contain in the closure the unique orbit of pencils
with two points of multiplicity $(k+1)$. \qed
\end{corollary0}

\section{The action on $MI_{\PP^3}(k)$ and $MI^s_{\PP^3}(k)$}

Along this section, we will keep the notations  introduced
 before. Since a morphism between $k$-instanton bundles lifts to
a morphism
between the corresponding monads, we get that the action of
$GL(V)$ over $MI_{\PP^3}(k)$ lifts to
the action of $GL(I)\times Sp(W)\times GL(V)$ over
${\cal Q}^0$. This means that if $E$ and $E'$ are $k$-instanton
bundles in the same $GL(V)$-orbit then any two representatives $A$,
$A'$ $\in{\cal Q}^0$ are in the same $GL(I)\times Sp(W)\times GL(V)$-orbit.
The (connected components containing the identity of the) stabilizers and the
stable points of the two actions correspond
to each other (this technique was used in \cite{AO99}).

As a consequence of the main theorem in \cite{CO00} we get

\begin{thm0}
All points in $MI_{\PP^3}(k)$ are semistable for the $GL(V)$-action.
\end{thm0}
\noindent {\bf Proof.}  Since ${\cal Q}^0$ is
the complement in ${\cal Q}$ of an invariant hypersurface
(see \cite{CO00}) all points in ${\cal Q}^0$ are semistable. In addition,
the two GIT-quotient $MI_{\PP^3}(k)/GL(V)$ and
${\cal Q}^0/GL(I)\times Sp(W)\times GL(V)$ are isomorphic. \qed

\begin{remark0}
\rm
\label{kernel} Let $U$ and $U'$ be vector spaces of dimension $2$.
Fix three isomorphisms $V\simeq U\otimes U'$,
$W\simeq S^{k}U\otimes U'$ and $I\simeq S^{k-1}U$. Then, there is an induced
  natural multiplication map
$A\in Hom(W\otimes\OO,I\otimes\OO(1))$ which is
$GL(U)\times GL(U')$-invariant. Indeed it can be shown that
$\ker A$ is a rank $k+2$ bundle on $\PP(V)$ with symmetry group in $GL(V)$
isomorphic to  $GL(U)\times GL(U')/\CC^*$. In the case $U=U'$,
 this construction was considered in \cite{ST90}, where $\ker A$
was called the kernel bundle. We could work in $SL(V)$, which is the
universal cover of $Aut(\PP(V))$, and in this case the symmetry group of
$\ker A$ is $Spin(V)\simeq SL(U)\times SL(U')$, but we prefer
to consider the $GL$ groups, which is a not essential variation.
\end{remark0}

After this remark, we are led to pose the following definition, which will
play an important role in the sequel.

\begin{defn0}
\label{def:properly} A subgroup $G\subseteq GL(V)$ is said to have
a fine action over $MI_{\PP^3}(k)$ if there are vector spaces $U$
and $U'$ of dimension $2$ such that $G\subseteq GL(U)\times
GL(U')$ and there exists a lifted action of $G$ over $W\otimes
I\otimes V$ which is determined by some fixed isomorphisms
$V\simeq U\otimes U'$, $W\simeq S^{k}U\otimes U'$, $I\simeq
S^{k-1}U$.
\end{defn0}

\begin{thm0}
\label{zkunst}
Let $E$ be a $k$-instanton bundle such that $Z_k(E)\neq\emptyset$.
Then $[E]$ is a not stable point in $MI_{\PP^3}(k)$ for the $GL(V)$-action.
\end{thm0}
\noindent {\bf Proof.} By \cite{Ski97}, $E$ is represented by
\[A=\left[I(x_0,x_1)|H(x_0,x_1,x_2,x_3)\right]\]
where $H$ is a $k \times (k+1)$ matrix of linear forms.
Let $\lambda\simeq\CC^*\subseteq SL(U')$ acting in such a way that
the splitting $W={U'}^{k+1}$ defines two eigenspaces in
$W$ with positive and negative weights that correspond to the two
blocks in which
are divided the columns of $A$. Moreover, the positive eigenspace on
 $V=U'\oplus U'$ is generated by $x_0, x_1$ and the negative one by
  $x_2, x_3$. The action on $I$ is trivial. $\CC^*$ has a fine  action
   according to \defref{def:properly} and we get
\[\lim_{t\to 0}\lambda (t)\cdot A=\left[I(x_0,x_1)|H(0,0,x_2,x_3)\right]\]
which still represents a $k$-instanton. \qed

\begin{remark0}
\rm By \cite{Rao97} we know that if  $Z_k(E)\neq\emptyset$ then
$[E]$ is a smooth point in $MI_{\PP^3}(k)$. It should be
interesting to know if the unstable points for the $GL(V)$-action
over $MI_{\PP^3}(k)$ are smooth points.
\end{remark0}

\vspace{3mm}

Let us see, by means of an example, that the converse of \thmref{zkunst}
 is no longer true.

\begin{example0}
\label{pitag} \rm
 Let
\[A=\begin{pmatrix}
  3x_{0} & x_{1} &  & 4x_{3} & &&5x_3 & x_{2} \\
  & 5x_{0} & x_{1} & && 5x_3& x_{2} &  \\
   && 5x_{0} & x_{1} & 3x_{3} & x_{2}&&4x_0\\
\end{pmatrix}.\]
It is easy to check that $AJA^t=0$ and that the associated morphism is
 surjective (the rank is $3$ on every $(x_0,\ldots, x_3)\in\PP^3$)
 so that $A$ represents
a $3$-instanton bundle $E$. Moreover $Z_3(E)=\emptyset$ and
$\CC^*\subseteq Sym (E)$ having a fine  action. The weights are as
follows:
  \begin{itemize}
  \item on $I: -2, 0, 2$
  \item on $W: 5,3,1,-1,-5,-3,-1,1$
  \item on $V: x_0 (-3), x_1(-1), x_2(1), x_3(3)$
\end{itemize}
which correspond to weights $-1$ and $1$ on $U$ and
$-2$ and $2$ on $U'$ with the isomorphisms
$I\simeq S^2U$, $W\simeq S^3U\otimes U'$, $V\simeq U\otimes U'$.
In particular $[E]$ is a not stable point for the $GL(V)$-action
over $MI_{\PP^3}(k)$ and this shows that the converse of \thmref{zkunst}
is not true.
\end{example0}

\begin{example0}
\label{ce} \rm
The generic $3$-instanton bundle $E$ with
$Z_3$ given by one point has $Sym^0(E)=0$.
This can be checked with the help of Macaulay system (\cite{BS})
by using Skiti monad (\cite{Ski97}) with a generic Hankel
matrix. Moreover,
it seems likely that the generic $k$-instanton bundle $E$ with
$Z_k\neq\emptyset$ has $Sym^0(E)=0$. On the other hand, there
are examples of $k$-instanton bundles such that
$dim~Sym(E)=2$.
\end{example0}

Keeping the notations introduced in Section 2,
 we recall  a result from \cite{CO00}.

\begin{thm0}
There is a $Sp(W)\times SL(I)\times SL(V)$-invariant homogeneous polynomial
$D$  over $W\otimes I\otimes V$ of degree $2k(k+1)$ such that
$A\in{\cal Q}^0$ if and only if $A\in{\cal Q}$ and $D(A)\neq 0$. \qed
\end{thm0}

\begin{corollary0}
\label{magic}
Let $A\in {\cal Q}^0$ representing a $k$-instanton bundle $E$ and
$\lambda\colon\CC^*\to Sp(W)\times SL(I)\times SL(V)$
be a morphism such that its image is contained in $Sym(E)$.
Then, if $\lim_{t\to 0}\lambda (t)\cdot A$ exists, it belongs
to ${\cal Q}^0$.
\end{corollary0}
\noindent {\bf Proof.} It follows from the fact that
$D(\lambda (t)\cdot A)$
is constant with respect to $t$. \qed

\begin{remark0}
\rm
After \corref{magic}, a classification of all $\CC^*$-invariant
$k$-instanton bundles
should give
the classification of unstable points in $MI_{\PP^3}(k)$, but we postpone this
study.
\end{remark0}

\vspace{3mm}

For special $k$-instanton bundles the description of the group action
is quite precise.
By the above geometric description, if $g^*E=E$, then $g$
leaves the smooth quadric $W(E)$ fixed and does not exchange the
two rulings. Hence, $g \in SL(U) \times SL(U')$, for some complex
vector spaces $U$ and $U'$ of dimension two, acting over
$\PP^{3}=\PP(U \otimes U')$ and one has to check how $ SL(U)
\times SL(U')$ acts on the space $G^k$.
Notice that the first $SL(U')$ does not change
anything so, for any  $ E \in MI_{\PP^3}^s(k)$
\refstepcounter{equation}
\[ (\theequation) \label{ii-iii} \hspace{7mm}
 SL(U') \subseteq Sym(E) \]
and it follows that
\refstepcounter{equation}
\[ (\theequation) \label{basiso} \hspace{7mm}
MI_{\PP^3}^s(k)/SL(V)\cong G^k/SL(U) \]
(compare it with \corref{isogk}). This description was performed for $k=2$
by Hartshorne
(\cite{Har78})  and in the general case by Spindler and Trautmann
(\cite{ST90}),
although they considered $U=U'$.

\vspace{3mm}

Let us now see that under this isomorphism, the isomorphism class of the bundle
is uniquely determined by the $SL(2)$-class of its associated Hankel matrix.
Indeed we have

\begin{lemma0}
\label{goodiso}
\begin{itemize}
\item[(i)]
 For every $g\in SL(2)$
\[ \left(\begin{array}{cc}(S^kg^{-1})^t&0\\
0&S^kg\end{array}\right)\in Sp(\CC^{2k+2})\]
with respect to the standard $J$.
\item[(ii)]
For any $(k+1) \times (k+1)$-Hankel matrix $H$
\[ \left(S^{k-1}g\right)^t\cdot
\left[ I_k(x_0, x_1) | I_k(x_2, x_3)\cdot H\right]\cdot
\left(\begin{array}{cc}(S^kg^{-1})^t&0\\
0&S^kg\end{array}\right)=
\left[ I_k(x'_0, x'_1) | I_k(x'_2, x'_3)\cdot H'\right] \]
being $H'=\left(S^kg\right)^t \cdot H \cdot S^kg$ a Hankel matrix.
\item[(iii)] The isomorphism $MI_{\PP^3}^s(k)/SL(V)\cong G^k/SL(U)$ takes
the $SL(V)$ class of the bundle to the
 $SL(2)$-class of its associated Hankel matrix given in (\ref{def:balanced}).
 \end{itemize}
\end{lemma0}

\noindent
{\bf Proof.} We will prove $(iii)$ since $(i)$ and $(ii)$ are easy to check.
Given $E$ a special $k$-instanton bundle
choose coordinates such that
$\{x_0x_3-x_1x_2=0\}\in W(E)$.
Take two $k$-jumping lines $L_1$ and $L_2$ such that
$H^0(E_{|L_1})$  and $H^0(E_{|L_2})$ are orthogonal spaces into $W$
with respect to $J$ (here $H^0(E_{|L_i})$ are considered into $W$ by the monad
(\ref{definst}) ).
Moreover we can assume that  $L_1=\{x_0=x_1=0\}$ and
$L_2=\{x_2=x_3=0\}$. With this choice of coordinates there exists
a matrix $A$ representing $E$ as in (\ref{def:balanced})
containing a Hankel matrix $H$.
$SL(2)$ acts on these choices, the $SL(2)$-class of $H$ is
uniquely determined and characterizes the $SL(V)$-orbit
of $E$.  \qed

\vspace{3mm}

As a consequence we obtain the following nice description

\begin{corollary0}
\label{semistable}
 Let $MI_{\PP^3}^{s}(k)$ be the moduli space of
special $k$-instanton bundles on $\PP^3$. Then all points of
$MI_{\PP^3}^{s}(k)$   are semistable for the action of $SL(V)$
and the only not stable ones have a section vanishing on a line
with multiplicity $(k+1)$ . All the not stable orbits contain in
the closure the unique orbit of bundles having two distinct
sections each one vanishing on a (different) line with
multiplicity $(k+1)$.
\end{corollary0}
\noindent {\bf Proof.} It follows from the isomorphism
(\ref{basiso}), \corref{isogk} and \lemref{goodiso}.
\qed

\begin{prop0}
\label{iii-ii} Let $E$ be a $k$-instanton bundle such that
$SL(U')\subseteq Sym(E)$ having a fine  action. Then $E$ is
special.
\end{prop0}
\noindent {\bf Proof.} By the assumption there is a
$\CC^*\subseteq SL(U')$ with a two-dimensional eigenspace in $V$
of positive weight. Then there is a matrix representing $E$
such that in the first $k\times (k+1)$ submatrix only the coordinates
of this eigenspace appear. It follows that the line spanned by this
coordinates is a $k$-jumping line. Since $Z_k(E)$ is $SL(U')$-invariant
and not empty, it is one-dimensional and the result follows from
\corref{line:special}. \qed


\vspace{3mm}

Let us now briefly describe different group actions of $G \subset Sym(E)$,
being $E$ a special $k$-instanton bundle.

Let $E$ be a special $k$-instanton bundle represented by a matrix $A$
(see \notref{defbala}). In general we have $SL(2)\subseteq Sym(E)$
acting in the following way.
If  \[    g= \left[  \begin{array}{cc}
     \alpha & \beta   \\
     \gamma & \delta \\
      \end{array} \right]\in SL(2) \]
then the action on $W$ is determined by
\[    \left[  \begin{array}{cc}
     \alpha\cdot Id& \beta \cdot H  \\
     \gamma\cdot H^{-1} & \delta \cdot Id\\
      \end{array} \right]\in Sp(W) \]
while the action on $I$ is trivial and the action on $V$ is given by
\[      \begin{array}{ccc}
x_0 &\mapsto& \alpha x_0+\gamma x_2   \\
x_1 &\mapsto& \alpha x_1+\gamma x_3   \\
x_2 &\mapsto& \beta x_0+\delta x_2   \\
x_3 &\mapsto& \beta x_1+\delta x_3.   \\
\end{array}  \]
Let us now describe some particular case. Consider
the Hankel matrix $H$ given by
\[ \alpha_i=\delta_{i,k}. \]

It corresponds to the form $x^ky^k$ which, for $k\ge 2$, is the only form such
that $\Delta \neq 0$ and such that the stabilizer has dimension $4$.
In this case,
$\CC^*\cdot SL(2)\subseteq Sym(E)$
where $SL(2)$ acts as above and $\CC^*$
acts in the following way. For $t\in\CC^*$
we have
\[    \left[  \begin{array}{ccc}
     t^{k-1} &&  \\
     &\ddots & \\
&&t^{-k+1}
      \end{array} \right]\in SL(I) , \hspace{3mm}
    \left[  \begin{array}{cccccc}
     t^{-k} &&&&&  \\
     &\ddots &&&&\\
&&t^{k}&&&\\
&&&t^k&&\\
&&&&\ddots&\\
&&&&&t^{-k}\\
      \end{array} \right]\in Sp(W) \]
and the action on $V$ is determined by
\[      \begin{array}{ccc}
x_0 &\mapsto& t x_0  \\
x_1 &\mapsto& t^{-1}x_1   \\
x_2 &\mapsto& t x_2   \\
x_3 &\mapsto& t^{-1} x_3.   \\
\end{array}  \]

A matrix description of the instanton bundle $E$ with
$Sym(E)=\CC^*\cdot SL(2)$ is
\refstepcounter{equation}
\[ (\theequation) \label{matriu} \hspace{7mm}
 \left ( \begin{array}{cccccccccc}
 x_0 & x_1 &       &&        &&&x_3&x_2\\
     & .& .   &  &&&.&.     \\
 &&.&.  &  &.&.&     \\
   &  &     & x_0&x_1 &x_3&x_2\\
     \end{array} \right ). \]

\vspace{4mm}

Finally, let us mention that bundles $E$ having two distinct
sections each one vanishing on a (different) line with
multiplicity $(k+1)$ are quite interesting, because they are
characterized by the property $Sym(E)\simeq SL(U)\cdot\CC^*$ with
a fine  action.

\vspace{4mm}

We will end this section given a more precise description of the
correspondence,
introduced in Section 2,
between special $k$-instanton bundles and linear systems $g_{k+1}^1$ without
base points. More precisely we will prove

\newpage

\begin{thm0}
\label{special}
 Let $MI_{\PP^3}^{s}(k)$ be the moduli space of
special $k$-instanton bundles on $\PP^3$. There is a natural
morphism
\[ \phi: MI_{\PP^3}^{s}(k) \longrightarrow \PP^{2k}/SL(U) \]
which takes $E$ to the branch locus defined by its pencil of sections
and which factors through
\[ \phi': MI_{\PP^3}^{s}(k)/SL(V) \longrightarrow \PP^{2k}/SL(U). \]
$\phi'$ is finite of degree equal to
$deg(G(\PP^1,\PP(S^{k+1}U)))=\frac{1}{k+1}\binom{2k}{k}=c_{k+1}$,
where $ c_{k+1}$ is the $(k+1)$-th Catalan number (see e.g. \cite{GKZ}; pag.
239).
\end{thm0}

\vspace{4mm}

The above theorem  appears in \cite{New81} for $k=2$ and
in this case $MI_{\PP^3}(2)/SL(V)$ is isomorphic to  $\AAA^1$.

 This result answers the question posed by R. Hartshorne in \cite{Har78b}
concerning the relation between the cross-ratio of the four branch
points of a $g_{3}^{1}$ and the orbits for the action of $SL(U)$.
Indeed, Hartshorne asked if the cross-ratio determined uniquely
the orbit and P. Newstead showed that for each value of the
cross-ratio there are exactly two orbits in $G^2$ with one
exception (See \cite{New81} for more details).

\vspace{4mm}
\noindent {\bf Proof of \thmref{special}.}

\vspace{4mm}

By the isomorphism (\ref{basiso}), $MI_{\PP^3}^{s}(k)/SL(V) \cong
G^k/SL(U)$. Moreover, by \lemref{goodiso}, this isomorphism takes
the $SL(V)$-class of the bundle to the $SL(2)$-class of its
associated Hankel matrix given in (\ref{def:balanced}). Hence, we
will prove that there exists a finite morphism
\[ \begin{array}{ccc}
\psi: G^k/SL(U) & \longrightarrow & \PP^{2k}/SL(U) \\
\end{array} \]
of degree equal to $c_{k+1}$. To this end, let $G=Gr(\PP^1,
\PP(S^{k+1}U))$ and denote by $p_{i,j}$ the Pl\"{u}cker
coordinates on $G$. Define
\[ \begin{array}{cccc}
R: &  G & \longrightarrow & \PP^{2k}=\PP(S^{2k}U) \\  & (p_{0,1},
\cdots, p_{k,k+1}) & \longrightarrow & (q_1: \cdots : q_{2k+1}) \\
\end{array} \]
where, for any $1 \leq m \leq 2k+1$, $q_m= \sum_{\mu+\nu=m,
\hspace{1mm} \mu < \nu} (\nu-\mu)p_{\mu,\nu}$ (compare it with
\cite{Tra88}; Pag. 41-42).

\noindent {\em Claim:} $R$ is induced by the projection from the
linear subspace defined by
\[ q_0= \cdots = q_{2k+1}=0\]
which is disjoint from $G$.

\noindent {\em Proof of the Claim:} First of all recall that, for
any $\{i_0,i_1,i_2,i_3 \} \subset \{1, \cdots, k\}$,  the
Pl\"{u}cker coordinates $p_{i,j}$ verify the following Pl\"{u}cker
relations
\[p_{i_0,i_1}p_{i_2,i_3}-p_{i_0,i_2}p_{i_1,i_3}+p_{i_0,i_3}p_{i_1,i_2}=0.\]

We will prove by induction on $m=\mu+\nu$, that if $q_1= \cdots =
q_{2k+1}=0$ then $p_{i,j}=0$ for any pair $(i,j)$. If $m=1$, the
assumption $q_1=0$ implies $p_{0,1}=0$. Let us assume
$p_{\mu,\nu}=0$ for $\mu < \nu$ such that $\mu + \nu=m$ and we
will see that $p_{\mu,\nu}=0$ for any $\mu < \nu$ such that $\mu +
\nu=m+1$. Considering the Pl\"{u}cker relations for a suitable
sets of indices and by induction hypothesis we get
\[ p_{i,j}p_{k,l}=0 \]
for any $i <j$, $k <l$ with $i+j=k+l=m+1$. Hence, from the
assumption
\[ 0=q_{m+1}= \sum_{\mu+\nu=m+1, \hspace{1mm}  \mu <
\nu} (\nu-\mu)p_{\mu,\nu} \] we deduce that $p_{\mu,\nu}=0$ for
any $\mu < \nu$ with $\mu+\nu=m+1$, which proves what we want.

It follows from the claim that $R$ is a finite morphism. Moreover,
since $G^k \subset G$, by \cite{Tra88}; Pag. 41, $R$ induces a
finite morphism
\[ \begin{array}{ccc}
\psi: G^k/SL(U) & \longrightarrow & \PP^{2k}/SL(U) \\
\end{array} \]
of degree equal to $c_{k+1}$ which, by \thmref{trautmann}, sends
the $SL(U)$-class of the linear system $g_{k+1}^1$ associated to
the $SL(V)$-class of a bundle $E$, to the branch locus defined by
its pencil of sections. \qed

\vspace{3mm}


\begin{remark0} \rm
A general special instanton bundle with $c_2=2$ has symmetry group
$SL(2)\cdot G_8$ where $G_8$ is binary dihedral of order $8$.
There is a central extension
\[ 0\rightarrow\ZZ_2\rightarrow G_8\rightarrow\ZZ_2\oplus\ZZ_2\rightarrow 0. \]
In fact $\ZZ_2\oplus\ZZ_2\cong D_4$. If the quartic form
corresponding to $(\alpha_0,\ldots ,\alpha_4)$ satisfies $I=0$
then $Sym(E)\cong  SL(2)\cdot G_{24}$ where $G_{24}$ is binary
tetrahedral of order $24$. There is a central extension
\[ 0\rightarrow\ZZ_2\rightarrow G_{24}\rightarrow A_4\rightarrow 0. \]
If the bundle has only one section vanishing on a line with multiplicity $4$
then the quartic form has a double root and
$Sym(E)\cong  SL(2)\cdot G_{16}$ where
$G_{16}$ is binary dihedral of order $16$.
The bundle with two sections each one vanishing on a different line
with multiplicity $4$ has
 $Sym(E)\cong  SL(2)\cdot \CC^*$ and it is characterized by the property
$\dim Sym(E)\ge 4$.
\end{remark0}

\section{Link between $Sym$ and $Z_k$}

The goal of this section is to study the symmetry group $Sym(E)$ of
a $k$-instanton bundle $E$ such that $Z_k(E) \neq \emptyset$.

\vspace{4mm}

\begin{thm0}
\label{big:case} Let $E$ be a $k$-instanton bundle on $\PP^{3}$
such that $Z_k(E)$ contains a double point. Then there exists $\CC
\subseteq Sym^0(E) $, where $\CC$ has a fine  action.
\end{thm0}
\noindent {\bf Proof.}
By assumption and \cite{Ski97} there is monad representing $E$ with
\[A=
\left( \left.\begin{array}{llll}
 x_0 & x_1 &       &        \\
     & \ddots& \ddots   &       \\
     &     & x_0&x_1 \\
     \end{array}\right| H(x_0, x_1, x_2, x_3) \right) \]
and the line $\{x_0=x_1=0\}\in Z_k(E)$ with
coordinates
$(p_{01},\ldots p_{13}, p_{23})=(0,\ldots ,0,1)$ corresponds to a
double point.
By  \remref{quadric}, $Z_k(E)$ is given by $\sum_{i<j}p_{ij}A_iJA_j^t=0$.
The crucial point is that
\[A'=
\left( \left.\begin{array}{llll}
 x_0 & x_1 &       &        \\
     & \ddots& \ddots   &       \\
     &     & x_0&x_1 \\
     \end{array}\right| H(0, 0, x_2, x_3) \right) \]
represents another instanton $E'$. Moreover
\[A_0JA_2^t={A'}_0J{A'}_2^t,\quad A_0JA_3^t={A'}_0J{A'}_3^t, \quad
 A_1JA_2^t={A'}_1J{A'}_2^t, \quad A_1JA_3^t={A'}_1J{A'}_3^t,\]
\[A_2JA_3^t={A'}_2J{A'}_3^t=0\]
and ${A'}_0J{A'}_1^t=0$ while, in general,  $A_0JA_1^t\neq 0$.
By assumption the variety $\sum_{i<j}p_{ij}A_iJA_j^t=0$ contains
a line tangent to the Pl\"ucker quadric
$p_{01}p_{23}-p_{02}p_{13}+p_{03}p_{12}=0$
in the point $(0,\ldots ,0,1)$.
Hence the system
\[p_{02}A_0JA_2^t+p_{03}A_0JA_3^t+p_{12}A_1JA_2^t+p_{13}A_1JA_3^t=0\]
has a nonzero solution $( \tilde{p}_{02}, \tilde{p}_{03},\tilde{
p}_{12},\tilde{p}_{13})$
which gives the line
$(0,s\tilde{p}_{02},s\tilde{p}_{03},s\tilde{p}_{12},
s\tilde{p}_{13},t)$ where $(s,t)\in\PP^1$.
Notice that this system is the same which gives solutions in
the unknowns $(p_{01},\ldots p_{13}, p_{23})$ for $Z_k(E')$
which has now the solutions
$(u,s\tilde{p}_{02},s\tilde{p}_{03},s\tilde{p}_{12},
s\tilde{p}_{13},t)$ where $(u,s,t)\in\PP^2$.
Therefore, $Z_k(E')$ is obtained by cutting the Pl\"ucker
 quadric with a plane which, in particular, means that it is a conic.
 Hence, by \corref{line:special}, $E'$ is special and
by the \lemref{balanced} there exists a nondegenerate
 $(k+1)\times (k+1)$-Hankel matrix $K$
such that
$I_k(x_2,x_3)\cdot K=H(0,0,x_2,x_3)$.
Define the morphism
\[\begin{array}{ccc}\CC&\to &Sp(W)\\
t&\mapsto & S_t\\
\end{array}\]
where $S_t$ is represented by the matrix
\[S_t=\left[\begin{array}{cc}Id& tK\\
0&Id\\ \end{array}\right].\]
We have
\[A'\cdot S_t=
\left( \left.\begin{array}{llll}
 x_0 & x_1 &       &        \\
     & \ddots& \ddots   &       \\
     &     & x_0&x_1 \\
     \end{array}\right| H(0, 0, tx_0+x_2, tx_1+x_3) \right). \]
Hence
\[A\cdot S_t=
\left( \left.\begin{array}{llll}
 x_0 & x_1 &       &        \\
     & \ddots& \ddots   &       \\
     &     & x_0&x_1 \\
     \end{array}\right| H(x_0, x_1, tx_0+x_2, tx_1+x_3) \right) \]
and
$\CC\subseteq Sym^0(E)$ as we wanted. \qed

\vspace{3mm}

\begin{remark0}
\rm
With the notations of the above proof, there is a subgroup
$\lambda$ isomorphic to $\CC^*$ given by $x_0\mapsto tx_0$, $x_1\mapsto tx_1$,
$x_2\mapsto t^{-1}x_2$, $x_3\mapsto t^{-1}x_3$
such that, by  \corref{magic}, $\lim_{t\to 0}\lambda (t)\cdot E$
is still a $k$-instanton bundle  with four-dimensional symmetry group as in
(\ref{matriu}).
\end{remark0}

\vspace{3mm}

\begin{thm0}
\label{twopoints} Let $E$ be a $k$-instanton bundle on $\PP^{3}$
such that $Z_k(E)$ contains two distinct points. Then there exists
$\CC^* \subseteq Sym^0(E) $, where $\CC^*$ has a fine  action.
\end{thm0}
 \noindent {\bf Proof.}
By assumption and \cite{Ski97}, there is monad representing $E$ with
\[A=
\left( \left.\begin{array}{llll}
 x_0 & x_1 &       &        \\
     & \ddots& \ddots   &       \\
     &     & x_0&x_1 \\
     \end{array}\right| H(x_0, x_1, x_2, x_3) \right) \]
and the line $L_1=\{x_0=x_1=0\}\in Z_k(E)$. Since $Z_k(E)$ contains two
distinct points, we can assume that $L_2=\{x_2=x_3=0\}\in Z_k(E)$.
Moreover, from the fact that $h^0(E|L_2)=k+1$ we deduce that there exists a
nondegenerate $(k+1) \times (k+1)$-Hankel matrix $M$ such that
\[  H(x_0,x_1,0,0)=I_k(x_0,x_1) \cdot M. \]
Defining
\[S=\left[\begin{array}{cc}Id& -M\\
0&Id\\ \end{array}\right]\]
we get
\[A \cdot S= \left [ I_k(x_0,x_1) | -M \cdot I_k(x_0,x_1)+H(x_0,x_1,x_2,x_3)
\right ]=  \left [ I_k(x_0,x_1) | \tilde{H}(x_2,x_3) \right ]\]
where $\tilde{H}(x_2,x_3)$ is a matrix of linear forms only in
$x_2,x_3$. Therefore, there exists $\CC^*$ acting on $U$ with
weights $-1,0$ and on $U'$ with weights $0,0$ such that $\CC^*
\subset Sym(E)$ having a fine  action, which proves what we want.
\qed

 \vspace{3mm}

\label{biblio}

\begin{tabular}[t]{lcl}

Laura Costa && Giorgio Ottaviani\\
Dept. Algebra y Geometr{\'\i}a & & Dipartimento di Matematica U.
 Dini \\
 Universitat de Barcelona & & Universita di Firenze \\
 Gran Via,585 & & Viale Morgagni, 67 A \\
  08007 Barcelona & & 50134 Firenze \\
  Spain & & Italy \\
 e-mail: costa@mat.ub.es  & &  e-mail: ottavian@math.unifi.it
\end{tabular}

\end{document}